%% file: 0-cdc21.tex
\documentclass[letterpaper, 10 pt, conference]{ieeeconf}  

\IEEEoverridecommandlockouts                              
\usepackage{graphicx} 
\usepackage{amsmath} 
\usepackage{amssymb}  
\usepackage{color,soul}
\usepackage[dvipsnames]{xcolor}
\usepackage{algorithm}
\usepackage{mathtools}
\usepackage{algpseudocode}
\algtext*{EndFor}
\algtext*{EndIf}

\usepackage{subcaption} 
\usepackage{hyperref}
\usepackage{enumerate}
\let\labelindent\relax
\usepackage[shortlabels]{enumitem}
\usepackage{cite}

\usepackage{tikz}
\usetikzlibrary{calc, shapes, backgrounds}
\usepackage{pgfpages}
\DeclareMathAlphabet{\mathcal}{OMS}{cmsy}{m}{n}
\usepackage{pgf}
\usetikzlibrary{arrows,automata,calc,positioning}

\usepackage{xcolor}
\definecolor{myLightBlue}{rgb}{0.357,0.627,0.796}
\definecolor{myMidBlue}{rgb}{0,0.4,0.663}
\definecolor{matlabBlue}{rgb}{0,0.447,0.741}
\definecolor{myRed}{rgb}{0.918,0.133,0.047}
\definecolor{myGrey}{rgb}{0.75,0.75,0.75}
\definecolor{myLightGrey}{rgb}{0.9,0.9,0.9}
\definecolor{myOrange}{rgb}{0.85,0.325,0.098}
\definecolor{myGreen}{rgb}{0.466,0.674,0.188}
\definecolor{myPurple}{rgb}{0.5,0.3,1}
\definecolor{mygreen}{RGB}{204,255,153}
\definecolor{myblue}{RGB}{153,204,255}
\definecolor{myorange}{RGB}{255,204,153}
\definecolor{mylightgray}{RGB}{225,225,225}

\newtheorem{defn}{\bf Definition}
\newtheorem{prop}{\bf Proposition}[section]
\newtheorem{lemma}{\bf Lemma}[section]
\newtheorem{example}{\bf Example}[section]

\newtheorem{thm}{\bf Theorem}[section]
\newtheorem{corol}{\bf Corollary}[section]
\newtheorem{lemmaappendix}{\bf Lemma}[subsection]

\providecommand{\com[2]}{\begin{tt}[#1: #2]\end{tt}}

\newcommand{\new}[1]{#1}

\newcommand{\newg}[1]{#1}

\newcommand{\newo}[1]{#1}

\providecommand{\prt}[1]{\ensuremath{\left( #1 \right)}}
\newcommand{\cqfd}{\hfill \rule{2mm}{2mm}\medbreak\indent}
\newcommand{\evec}{\mathbf{e}}
\DeclareMathOperator{\newrank}{rank} 
\newcommand{\rank}{\newrank\,}
\DeclareMathOperator{\newgrank}{generic\ rank} 
\newcommand{\grank}{\newgrank\,}
\DeclareMathOperator{\newvec}{vec} 
\newcommand{\myvec}{\newvec}

\DeclareMathOperator{\newsgn}{sgn} 
\newcommand{\sgn}{\newsgn}

\title{\LARGE \bf
Path-Based Conditions for Local Network Identifiability -- Full Version
}

\author{Antoine Legat and Julien M. Hendrickx
\thanks{A. Legat and J. M. Hendrickx are with ICTEAM Institute, UCLouvain, Belgium.
Work supported
by the ``RevealFlight'' ARC at UCLouvain, and by the Incentive Grant
for Scientific Research (MIS) ``Learning from Pairwise Data'' of the
F.R.S.-FNRS.{\tt\small antoine.legat@uclouvain.be, julien.hendrickx@uclouvain.be.}}%
}

\begin{document}

\maketitle
\thispagestyle{empty}
\pagestyle{empty}

\begin{abstract}
This work focuses on the generic identifiability of dynamical networks with partial excitation and measurement: a set of nodes are interconnected by transfer functions according to a known topology, some nodes are excited, some are measured, and only a part of the transfer functions are known. Our goal is to determine whether the unknown transfer functions can be generically recovered based on the input-output data collected from the excited and measured nodes.

We propose a decoupled version of generic identifiability that is necessary for generic local identifiability and might be equivalent as no counter-example to sufficiency has been found yet in systematic trials. This new notion can be interpreted as the generic identifiability of a larger network, obtained by duplicating the graph, exciting one copy and measuring the other copy. We establish a necessary condition for decoupled identifiability in terms of vertex-disjoint paths in the larger graph, and a sufficient one.
\end{abstract}

\section{INTRODUCTION}
\input{1-intro}
\medskip

\section{LOCAL IDENTIFIABILITY}
\input{2-local-identif}
\medskip

\section{DECOUPLED IDENTIFIABILITY}
\label{sec:decoup-identif}
\input{3-decoup-identif}

\section{PATH-BASED CONDITIONS\\FOR DECOUPLED IDENTIFIABILITY}
\label{sec:conditions}
\input{4-conditions}

\section{DISCUSSION}
\label{sec:discussion}
\input{5-dicussion}
\medskip

\section{CONCLUSION}
\input{6-conclusion}
\medskip

\bibliographystyle{ieeetr}
\bibliography{0-cdc21}

\appendix
\section{APPENDIX}
\input{7-appendix}

\end{document}

%% file: 1-intro.tex
\medskip

This paper addresses the identifiability of dynamical networks in which node signals are connected by causal linear time-invariant transfer functions, and can be excited and/or measured. Such networks can be modeled as directed graphs where each edge carries a transfer function, and known excitations and measurements are applied at certain nodes.

\medskip

We consider the identifiability of a network matrix $G(q)$, where the network is made up of $n$ node signals $w(t)$, external excitation signals $r(t)$, measured nodes $y(t)$ and noise $v_1(t), v_2(t)$ related to each other by:
\begin{align} \label{eq:networkModel} \begin{split}
	w(t) &= G(q)w(t) + Br(t) + v_1(t)\\
	y(t) &= Cw(t) + v_2(t),
\end{split} \end{align}
where matrices $B$ and $C$ are binary selections defining respectively the $n_B$ excited nodes  and $n_C$ measured nodes, forming sets $\mathcal{B}$ and $\mathcal{C}$ respectively. Matrix $B$ is full column rank and each column contains one 1 and $n-1$ zeros. Matrix $C$ is full row rank and each row contains one 1 and $n-1$ zeros.
The nonzero entries of the network matrix $G(q)$ define the network topology: some of them are known and collected in $G^0(q)$, and the others are the unknowns to identify, collected in $G^\Delta(q)$, such that $G(q) = G^0(q) + G^\Delta(q)$. The known edges are collected in set $E^0$, the unknown ones in $E^\Delta$, and $E = \{E^0,  E^\Delta\}$ is the set of all edges.

\medskip
		
We assume that \emph{the input-output relations between the excitations $r$ and measurements $y$ have been identified}, and that the network topology is known.
From this knowledge, we aim at recovering the unknown entries of $G(q)$.

\medskip

The model \eqref{eq:networkModel} has recently been the object of a significant research effort.
If the whole network is to be recovered, the notion of network identifiability is used \cite{weerts2015identifiability}. If one is interested in identifying a single module, topological conditions are derived in \cite{weerts2018single, gevers2018practical}. In this paper, we do not consider the impact of noise $v_1,v_2$, but studying the influence of rank-reduced or correlated noise under certain assumptions yields less conservative identifiability conditions \cite{weerts2018prediction, gevers2018identifiability, van2019local, ramaswamy2020local}.

\medskip

An approach dual to ours is to assume that network dynamics are known, and aim at identifying the topology from input/output data. This problem is referred to as topology identification, and is addressed in e.g. \cite{van2019topology_heterogeneous, van2019topology_lyapunov}.
	
\medskip
	
It turns out that the identifiability of the network, i.e. the ability to recover a module or the whole network from the input-output relation, is a generic notion: Either \emph{almost all} transfer matrices corresponding to a given network structure are identifiable, in which case the structure is called \emph{generically identifiable}, or none of them are. 
A number of works study generic identifiability when all nodes are excited or measured, i.e. when $B$ or $C=I$ \cite{bazanella2017identifiability, weerts2018identifiability}. Considering the graph of the network, path-based conditions on the allocation of measurements/excitations in the case of full excitation/ measurement are derived in \cite{hendrickx2018identifiability} /\cite{shi2020excitation}. Reformulating these conditions by means of disjoint trees in the graph, \cite{cheng2019allocation} develops a scalable algorithm to allocate excitations/measurements in case of full measurement/excitation. In case of full measurement, \cite{shi2020generic} derives path-based conditions for the generic identifiability of a subset of modules, under the presence of noise. Abstractions of dynamic networks yield conditions on nodes to measure for identifiability of a target module \cite{weerts2020abstractions}.
	
\medskip
	
The conditions of \cite{hendrickx2018identifiability} apply for generic identifiability, i.e. identifiability of almost all transfer matrices corresponding to a given network structure.  \cite{van2018topological} extends the path-based conditions under full excitation to determine the identifiability \emph{for all} (nonzero) transfer matrices corresponding to a given structure, and \cite{van2019necessary} provides conditions for the outgoing edges of a node, and the whole network under the same conditions.
	
\medskip
	
In all these works, the common assumption is that all nodes are either excited, or measured. In \cite{bazanella2019network}, this assumption is relaxed and generic identifiability with partial excitation and measurement is addressed for particular network topologies.
Partial excitation and measurement is also addressed by \cite{shi2020single}, which derives conditions for generic identifiability in terms of disconnecting sets, akin to what was done in full excitation/measurement \cite{shi2020excitation}.
	
\medskip

In the general case of arbitrary topology, partial excitation and measurement, \cite{legat2020local} introduces the notion of local identifiability, i.e. only on a neighborhood of $G(q)$. Local identifiability is a generic property, necessary for generic identifiability and no counterexample to sufficiency is known to the authors, i.e. no network which is locally identifiable but not globally identifiable. \cite{legat2020local} derives algebraic necessary and sufficient conditions for generic local identifiability for both the whole network, and a single module.

\medskip

The algebraic conditions of \cite{legat2020local} allow rapidly testing local identifiability for any given network, but we wish to find a combinatorial characterization for generic identifiability, that is expressed purely in terms of graph-theoretical properties, akin to what was done in the full excitation case e.g. in \cite{hendrickx2018identifiability}. Such characterization would in particular pave the way for optimizing the selection of nodes to be excited and measured, akin to the work in \cite{shi2020excitation} in the full measurement case. \cite{shi2020single} already provides conditions in terms of disconnecting sets, but we believe that the vertex-disjoint paths conditions of \cite{hendrickx2018identifiability} can be extended to local identifiability under partial excitation and measurement.

\medskip

In this paper, we derive path-based local identifiability conditions from the results of \cite{legat2020local}, in the general case of arbitrary topology, partial excitation and measurement. We extend the results of \cite{legat2020local} when some transfer functions are known \emph{a priori}. A more general notion of local identifiability is introduced: \emph{decoupled} identifiability, necessary for local identifiability. Interestingly, no counterexample to sufficiency is known to the authors, despite extensive testing (code available at \cite{matlab}).
Then, necessary and sufficient path-based conditions for decoupled identifiability are derived. These conditions are given in terms of connected paths and vertex-disjoint paths, and extend what one had in the full excitation case e.g. in \cite{hendrickx2018identifiability}.

\medskip
\noindent\textbf{Assumptions}:
We consider the problem modeled in \eqref{eq:networkModel}. Consistently with previous works, we assume that the network is well-posed, \newg{that is $(I-G(q))^{-1}$ is proper and stable, }
and we assume that $CT(q)B=C(I-G(q))^{-1}B$ has been identified exactly.
We do not suppose having access to any information related to the effect of the noise signals $v_1,v_2$. The additional information that could be gathered from this knowledge in our context is left for further works.

\medskip

Consistently with \cite{legat2020local}, we consider in this paper a single frequency $z$, so that all transfer functions are modeled simply by a complex value, and the matrices $G$ and $T(G)= (I-G)^{-1}$ are complex matrices rather than matrices of transfer functions.
Conceptually, our generic results directly extend to the transfer function case: if one can recover a $G_{ij}(z)$ at a given frequency $z$ for almost all $G$ consistent with a network, then one can also recover it at all other frequencies, and hence recover the transfer function.
We intend to remove this simplification or to formalize this intuitive argument in a further version of this work. 
In the remainder of this document, we omit $(q)$ to lighten notations. Also, the proofs of this paper are collected in the Appendix.

%% file: 2-local-identif.tex
We start from the definition of identifiability, see e.g. {\cite{hendrickx2018identifiability}}, which we extend to the case where some transfer functions are known ($G^0$), and some are not ($G^\Delta$), as in \cite{weerts2018single}. In the remainder of this paper, we denote $T(G) = (I-G)^{-1}$ and we sometimes drop the $(G)$ when there is no ambiguity.

\medskip

\begin{defn} \label{def:global_identif}
	The transfer function $G_{ij}$ is \emph{identifiable} at $G$ from excitations $\mathcal{B}$ and measurements $\mathcal{C}$ if, for all network matrix $\tilde{G}$ with same zero and known entries as $G$, there holds
	\begin{align} \label{eq:def_identif_edge}
		C \, T(\tilde{G}) \, B = C \, T(G) \, B \Rightarrow \tilde{G}_{ij} = G_{ij}.
	\end{align}
	The network matrix is identifiable at $G$ if each unknown transfer function $G_{ij}$ is identifiable at $G$, i.e. if the left-hand side of {\eqref{eq:def_identif_edge}} implies $\tilde G^\Delta = G^\Delta$.
\end{defn}

\medskip

We remind a notion of identifiability amenable to linear analysis: local identifiability, which corresponds to identifiability provided that $\tilde G$ is sufficiently close to $G$. Again, we extend the definition of \cite{legat2020local} to the more general case where some transfer functions are already known.

\medskip

\begin{defn} \label{def:local_identif}
	The transfer function $G_{ij}$ is \emph{locally identifiable} at $G$ from excitations $\mathcal{B}$ and measurements $\mathcal{C}$ if there exists $\epsilon > 0$ such that for any $\tilde{G}$ with same zero and known entries as $G$ satisfying $||\tilde{G}-G||<\epsilon$, there holds
	\begin{align} \label{eq:def_local_identif_edge}
	    C \, T(\tilde{G}) \, B = C \, T(G) \, B \Rightarrow \tilde{G}_{ij} =G_{ij}.
	\end{align}
	The network matrix is locally identifiable at $G$ if each unknown transfer function $G_{ij}$ is locally identifiable at $G$, i.e. if the left-hand side of \eqref{eq:def_local_identif_edge} implies $\tilde G^\Delta = G^\Delta$.
\end{defn}

\medskip

As stressed in \cite{legat2020local}, local identifiability is a necessary condition for identifiability. It is \emph{a priori} a weaker notion, yet no example of network locally identifiable but not globally identifiable is known to the authors.
Moreover, one can show that if a network is locally identifiable, then it can be recovered up to a discrete ambiguity, i.e. the set of unknown $G^\Delta$ corresponding to a measured $CTB$ would be discrete.

\medskip
	  
\noindent\textbf{Genericity:} Given a network and sets $\mathcal{B},$ $\mathcal{C}$ of excited and measured nodes, we say that an edge is \emph{generically} (locally) identifiable if it is (locally) identifiable at all $G$ consistent with the graph and known transfer functions, except possibly those lying on a lower-dimensional set \cite{davison1977connectability, dion2003generic} (i.e. a set of dimension lower than $|E^\Delta|$). In the remainder of this paper, we say that a property is \emph{generic}
if it either holds (i) for \emph{almost all} variables, i.e for all variables except possibly those lying on a lower-dimensional set, or (ii) for no variable.

\newg{For example, take a polynomial $p$. The \emph{nonzeroness} of $p(x)$ is a generic property of $x$: either (i) $p(x) \neq 0$ for all $x$ except its roots, or (ii) $p$ is the zero polynomial, which returns zero for all $x$.}

A handy consequence of this definition is the following: showing that a generic property holds for one variable $x$ implies that it holds for almost all $x$.

\subsection{Algebraic condition}

\medskip

Proposition \ref{prop:network_identif_rank_condition} below, adapted from \cite{legat2020local}, states that local identifiability is a generic property which can be checked by computing the rank of the matrix $K$:
\begin{align} \label{eq:def_K}
	K(G) \triangleq \prt{B^\top T^\top (G) \otimes C \, T(G)} I_{G^\Delta},
\end{align}
where $\otimes$ denotes the Kronecker product and the matrix $I_{G^\Delta} \in \{0,1\}^{n^2 \times |E|}$ selects only the columns of the preceding $n_B n_C \times n^2$ matrix corresponding to unknown edges:
\begin{align} \label{eq:IGDelta}
	I_{G^\Delta} \triangleq \left[ \myvec(G^\Delta(\evec_1)) \, \cdots \, \myvec(G^\Delta(\evec_{|E^\Delta|})) \right],
\end{align}
\newo{where $G^\Delta(\mathbf{v})$ is the matrix with same zero structure as $G^\Delta$, its nonzero entries collected in vector $\mathbf v$, }
 and $\evec_e$ is the standard basis vector filled with zeros except 1 at $e$-th entry.

\medskip

\begin{prop} \label{prop:network_identif_rank_condition}
Exactly one of the two following holds:\footnote{\newg{Observe that this implies: $K(G)$ has full rank for almost all $G$ \emph{if and only if} the network is generically locally identifiable, but the way Proposition {\ref{prop:network_identif_rank_condition}} is stated is stronger. The same phrasing remark holds for Proposition {\ref{prop:rank_hat_K_bilidentif}} and Corollary {\ref{corol:det_hat_K_bilidentif}}.}
}
	\begin{enumerate}[(i)]
	    \item $K(G)$ has full rank for almost all $G$ and the network is generically locally identifiable;
		\item $K(G)$ has full rank for no $G$ and the network is never locally identifiable.
	\end{enumerate}
	Moreover, $K(G)$ has full rank if and only if the following implication holds for all $\Delta$ with same zero entries as $G^\Delta$:
	\begin{align} \label{eq:CTDeltaTB_Delta_net}
		C \,  T(G)  \, \Delta  \, T(G) \,  B = 0 \Rightarrow \Delta = 0.
	\end{align}
\end{prop}

\medskip

The proof is given in \cite{legat2020local} when all transfer functions are unknown, and can straightforwardly be adapted when some transfer functions are already known.

\medskip

In this paper, we study network identifiability, i.e. the identifiability of \emph{all} unknown transfer functions $G^\Delta$.
Algebraic conditions for the identifiability of a single transfer function $G_{i,j}$ are derived in \cite{legat2020local}, but so far we were unable to interpret them in terms of paths in the graph. If the network is not identifiable, we do not have path-conditions to find out which transfer functions are problematic.

%% file: 3-decoup-identif.tex
\medskip

\begin{figure}
    \centering
    \includegraphics[width=\linewidth]{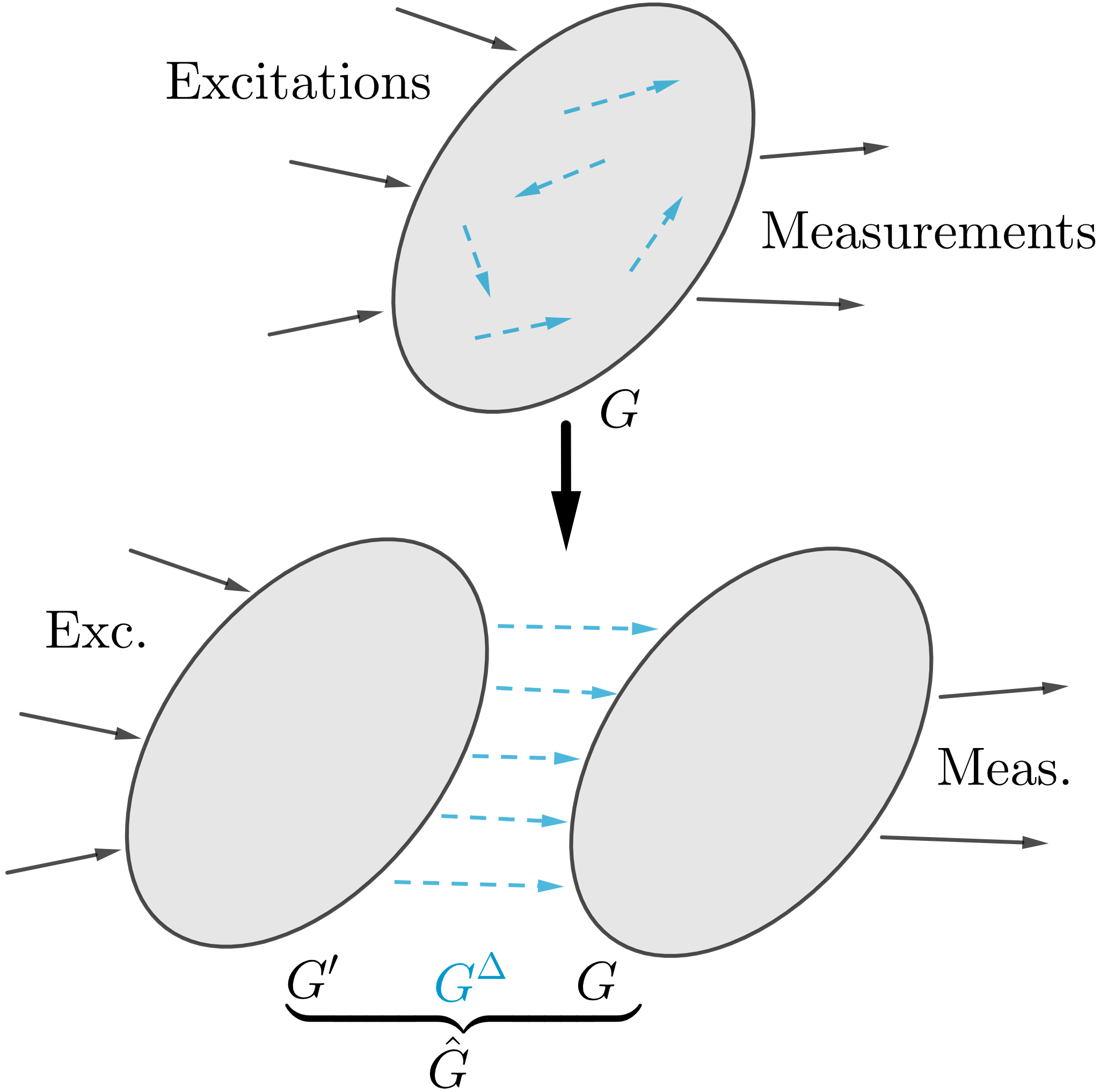}
    \caption{The decoupled network $\hat G$ related to $G$. The network is duplicated, excitations are applied on the left copy, measurements taken on the right, and unknowns in the middle.}
    \label{fig:decoupled_net}
\end{figure}

Consider condition \eqref{eq:CTDeltaTB_Delta_net} of Proposition \ref{prop:network_identif_rank_condition}.
We know that $C \, T(G) \, B$ is the closed-loop transfer matrix of the network with excitations $\mathcal{B}$ and measurements $\mathcal{C}$. It motivates the introduction of a larger network, whose closed-loop transfer matrix is given by $C \, T(G) \, \Delta \,  T(G) \, B$. It is built by duplicating the network, exciting the left copy, measuring the right copy and adding the unknown transfer functions in the middle (from left to right), see Fig. \ref{fig:decoupled_net}.

\medskip

Since the network is duplicated, we may allow the left and right copies to have different unknown transfer functions, while keeping the same topology and known transfer functions. We will see that relaxing the equality of left and right copies leads to a simpler analysis, hence we consider this more general notion of identifiability in this paper.

\newpage

\vspace*{-1cm}

First, we introduce the \emph{decoupled network}, whose closed-loop transfer matrix is $C \, T(G) \, \Delta  \, T(G') \, B$, see Fig. \ref{fig:decoupled_net}.




\medskip

\begin{defn} \label{def:decoupled_network}
    Consider a network of $n$ nodes with excitation matrix $B$, measurement matrix $C$ and network matrix $G = G^0 + G^\Delta$, where $G^0$ collects the known transfer functions and $G^\Delta$ collects the unknown transfer functions. Its \emph{decoupled network} is composed of $2n$ nodes: $\{1, \dots, n, 1', \dots, n'\}$. Its network matrix is defined by
    \begin{align*}
        \hat G (G,G') \triangleq
        \begin{bmatrix}
            G & G^\Delta \\
            0 & G'
        \end{bmatrix},
    \end{align*}
    where $G'$ has the same zero and known entries as $G$, \newo{and their unknown entries are given fixed parameters}. Transfer matrices $G$ and $G'$ are then fully known, while $G^\Delta$ contains the unknown transfer functions.
    Excitations are applied on the left subgraph ($G'$), and measurements on the right one ($G$), i.e. its excitation and measurement matrices are
    \begin{align*}
        \hat B \triangleq
        \begin{bmatrix}
            0 & 0 \\
            0 & B
        \end{bmatrix},
        && \hat C \triangleq
        \begin{bmatrix}
            C & 0 \\
            0 & 0
        \end{bmatrix}.
    \end{align*}
\end{defn}

\medskip 

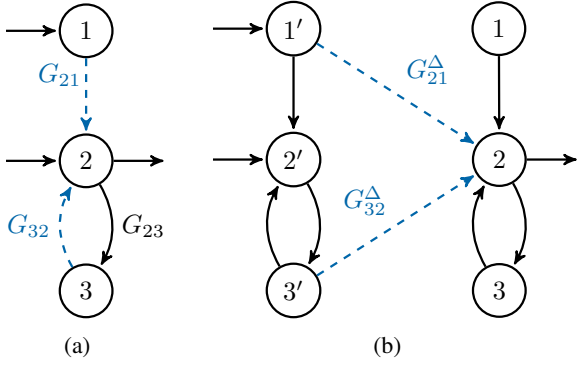
\begin{figure}
	\centering
	\subfloat[]{\label{subfig1}
	\begin{tikzpicture}[->,>=stealth',shorten >=1pt,auto,node distance=1.5cm,line width=0.3mm]
	  \node[state, minimum size=0.7cm] (N1)	                       {$1$};
	  \node[state, minimum size=0.7cm] (N2) [below = 1 cm of N1]   {$2$};
	  \node[state, minimum size=0.7cm] (N3) [below = 1 cm of N2]   {$3$};
	  \node        (N1e) [left = 0.7 cm of N1] {};
	  \node        (N2e) [left = 0.7 cm of N2] {};
	  \node        (N2m) [right = 0.7 cm of N2] {};
	  \path (N1) edge [above, dashed, color = myMidBlue] node {$G_{21}\qquad$} (N2)
			(N2) edge [bend left] node {$G_{23}$} (N3)
	 	   	(N3) edge [bend left, dashed, color = myMidBlue] node {$G_{32}$} (N2)
	 	   	(N1e) edge (N1)
	 	   	(N2e) edge (N2)
	 	   	(N2) edge (N2m);
	\end{tikzpicture}}
	\subfloat[]{\label{subfig2}
	\begin{tikzpicture}[->,>=stealth',shorten >=1pt,auto,node distance=3cm,line width=0.3mm]
	  \node[state, minimum size=0.7cm] (N1b)	            {$1'$};
	  \node[state, minimum size=0.7cm] (N2b) [below = 1 cm of N1b] {$2'$};
	  \node[state, minimum size=0.7cm] (N3b) [below = 1 cm of N2b] {$3'$};
	  \node[state, minimum size=0.7cm] (N1) [right = 2 cm of N1] {$1$};
	  \node[state, minimum size=0.7cm] (N2) [right = 2 cm of N2b] {$2$};
	  \node[state, minimum size=0.7cm] (N3) [right = 2 cm of N3b] {$3$};
	  \node        (N1be) [left = 0.7 cm of N1b] {};
	  \node        (N2be) [left = 0.7 cm of N2b] {};
	  \node        (N2m) [right = 0.7 cm of N2] {};
	  \path (N1b) edge (N2b)
			(N2b) edge [bend left] (N3b)
	 	   	(N3b) edge [bend left] (N2b)
	 	   	(N1) edge (N2)
			(N2) edge [bend left] (N3)
	 	   	(N3) edge [bend left] (N2)
	 	   	(N1b) edge [dashed, color = myMidBlue] node {$G_{21}^\Delta$} (N2)
	 	   	(N3b) edge [dashed, color = myMidBlue] node {$G_{32}^\Delta$} (N2)
	 	   	(N1be) edge (N1b)
	 	   	(N2be) edge (N2b)
	 	   	(N2) edge (N2m);
	\end{tikzpicture}}
	\caption{(a): A basic network example: unknown transfer functions are in dashed blue, arrows from the left denote excitations and the arrow to the right is a measurement. (b): Its decoupled network: The graph is duplicated, excitations are applied on left copy and measurement taken on the right one. The unknown transfer functions are in the middle, in dashed blue: they link the excited and measured subgraphs.}
	\label{fig:dummy_example_decoupled_T}
\end{figure}

An example of decoupled network is given in Fig. \ref{fig:dummy_example_decoupled_T}.

\medskip

We are now ready to introduce decoupled identifiability, which we will prove to be a generic property, necessary for generic (local) identifiability.
Take condition \eqref{eq:CTDeltaTB_Delta_net} of Proposition \ref{prop:network_identif_rank_condition}, and allow the two $T(G)$s to have different unknown transfer functions. The problem is no longer quadratic in $T(G)$, but linear in both $T(G')$ and $T(G)$:

\medskip

\begin{defn} \label{def:decoupled_identif}
    A network is \emph{decoupled-identifiable} at $(G,$ $G')$, with $G$ and $G'$ sharing the same zero and known entries, if for all $\Delta$ with same zero entries as $G^\Delta$, there holds:
    \begin{align} \label{eq:CTDeltaT'B_bili}
		C  \, T(G) \,  \Delta \,  T(G')  \, B = 0 \Rightarrow \Delta = 0.
	\end{align}
\end{defn}

\vspace*{-5mm}

Similarly to local identifiability, decoupled identifiability is a generic property: either it holds for almost all $(G,G')$, or for no $(G,G')$. It is proved in Proposition \ref{prop:rank_hat_K_bilidentif} below, which relies on the rank of the following matrix, defined analogously to $K$ in \eqref{eq:def_K}, but with $G'$ in the excitation part:
\begin{align*}
	\hat K(G,G') \triangleq \prt{B^\top T(G')^\top \otimes C \, T(G)} I_{G^\Delta},
\end{align*}
where $\otimes$ is the Kronecker product and $I_{G^\Delta}$ is defined in \eqref{eq:IGDelta}.

\medskip

\begin{prop} \label{prop:rank_hat_K_bilidentif}
Exactly one of the two following holds:
	\begin{enumerate}[(i)]
	    \item $\hat K(G,G')$ has full rank for almost all $(G,G')$ and the network is generically decoupled-identifiable;
		\item $\hat K(G,G')$ has full rank for no $(G,G')$ and the network is never decoupled-identifiable.
	\end{enumerate}
\end{prop}

\medskip

A first important proposition is that generic decoupled identifiability is necessary for generic local identifiability:

\medskip

\begin{prop} \label{prop:bili_necessary}
    If a network is generically locally identifiable, then it is generically decoupled-identifiable.
\end{prop}

\medskip

Interestingly, no counterexample to sufficiency is known to the authors, despite numerous tests:
we have randomly generated $10^6$ networks with up to $100$ nodes, and checked the generic local identifiability and generic decoupled identifiability of each network by computing $\rank K$ and $\rank \hat K$: for every network, both matched.
In other words, experiments seem to show that generic decoupled identifiability is equivalent to generic local identifiability.
Code available at \cite{matlab}.

\medskip

Moreover, since generic local identifiability is necessary for generic identifiability \cite{legat2020local}, so is generic decoupled identifiability. Hence, the necessary conditions derived for generic decoupled identifiability in this paper also hold for generic identifiability, and for identifiability. Besides, no example of network that is generically locally identifiable, but not generically identifiable is known to the authors.

\medskip

Definition \ref{def:decoupled_network} introduces the decoupled network, Definition \ref{def:decoupled_identif} presents the notion of decoupled identifiability. The proposition below unifies those two notions.

\medskip

\begin{prop} \label{prop:equivalence_bili_decoupled}
    The network $G$ is generically decoupled-identifiable if and only if its decoupled network $\hat G(G,G')$ is generically identifiable \newo{for almost all $(G,G')$.}
\end{prop}

\medskip

The following example illustrates propositions \ref{prop:rank_hat_K_bilidentif} and \ref{prop:equivalence_bili_decoupled}.

\medskip

\begin{example}
    For the network of Fig. \ref{fig:dummy_example_decoupled_T} (a), one can check that $\hat K$ has generic rank 2. Since there are 2 unknown transfer functions, Proposition \ref{prop:rank_hat_K_bilidentif} ensures that we have generic decoupled identifiability.
    It can be interpreted on the decoupled network, depicted on Fig. \ref{fig:dummy_example_decoupled_T} (b), as the generic identifiability of the unknown transfer functions $G^\Delta_{21}$ \& $G^\Delta_{32}$.
\end{example}

%% file: 4-conditions.tex
For ease of presentation, we consider the case where there are exactly $n_B n_C$ unknown edges, i.e. as many as the number of (excitation $b$, measurement $c$) pairs. If there are more unknown edges, then the network is not identifiable since there are more unknowns than (in, out) data.
The situation with less unknown edges than $n_B n_C$ will be addressed in a more complete version of this work.
In this section, we drop the $G,G'$ arguments and denote $T(G')$ by $T'$.
Also, we refer to unknown transfer functions $G^\Delta_{i,j}$ as unknown edges $\alpha$.

\medskip

Since $|E^\Delta| = n_B n_C$, $\hat K$ is a square matrix, hence $\rank \hat K = |E^\Delta|$ is equivalent to $\det \hat K \neq 0$. Proposition \ref{prop:rank_hat_K_bilidentif} can then be rewritten in terms of determinant:

\medskip

\begin{corol} \label{corol:det_hat_K_bilidentif}
Exactly one of the two following holds:
	\begin{enumerate}[(i)]
	    \item $\det \hat K \neq 0$ generically and $G$ is generically decoupled-identifiable;
		\item $\det \hat K$ is always zero and $G$ is never decoupled-identifiable.
	\end{enumerate}
\end{corol}

\medskip

In order to interpret $\det \hat K$, we develop the system of \eqref{eq:CTDeltaT'B_bili}:
\begin{align} \label{eq:developed_system}
    \sum_{\alpha \in E^\Delta} T'_{\alpha,b} T_{c,\alpha} \, \Delta_\alpha = 0 \qquad \forall \ b \in \mathcal{B}, \ c \in \mathcal{C}.
\end{align}
Each equation of \eqref{eq:developed_system} represents a pair (excitation $b$, measurement $c$), and corresponds to a row of $\hat K$. Each column of $\hat K$ matches an unknown edge $\alpha$, and the entry of $\hat K$ corresponding to $(b,c)$ and $\alpha$ is given by $T'_{\alpha,b} T_{c,\alpha}$.\footnote{\newo{By abuse of notation, $T'_{\alpha,b}$ stands for the transfer function between excitation $b$ and \emph{start node of} edge $\alpha$, and $T_{c,\alpha}$ denotes the one between \emph{end node of} edge $\alpha$ and measurement $c$.}
}

\medskip

\begin{figure}
    \centering
    \subfloat[]{\label{subfig1}\includegraphics[width=\linewidth]{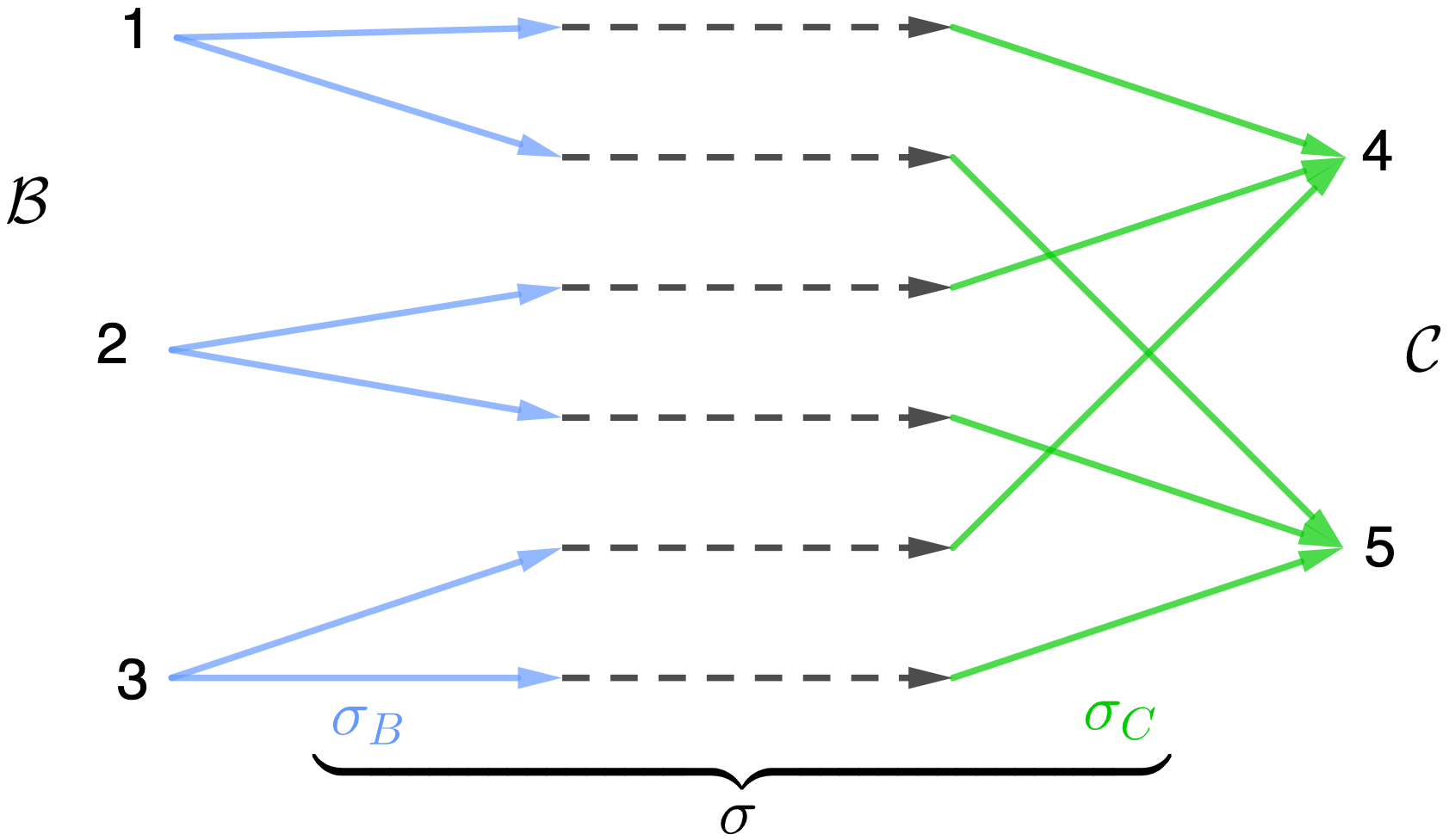}} \\
    \subfloat[]{\label{subfig2}\includegraphics[width=\linewidth]{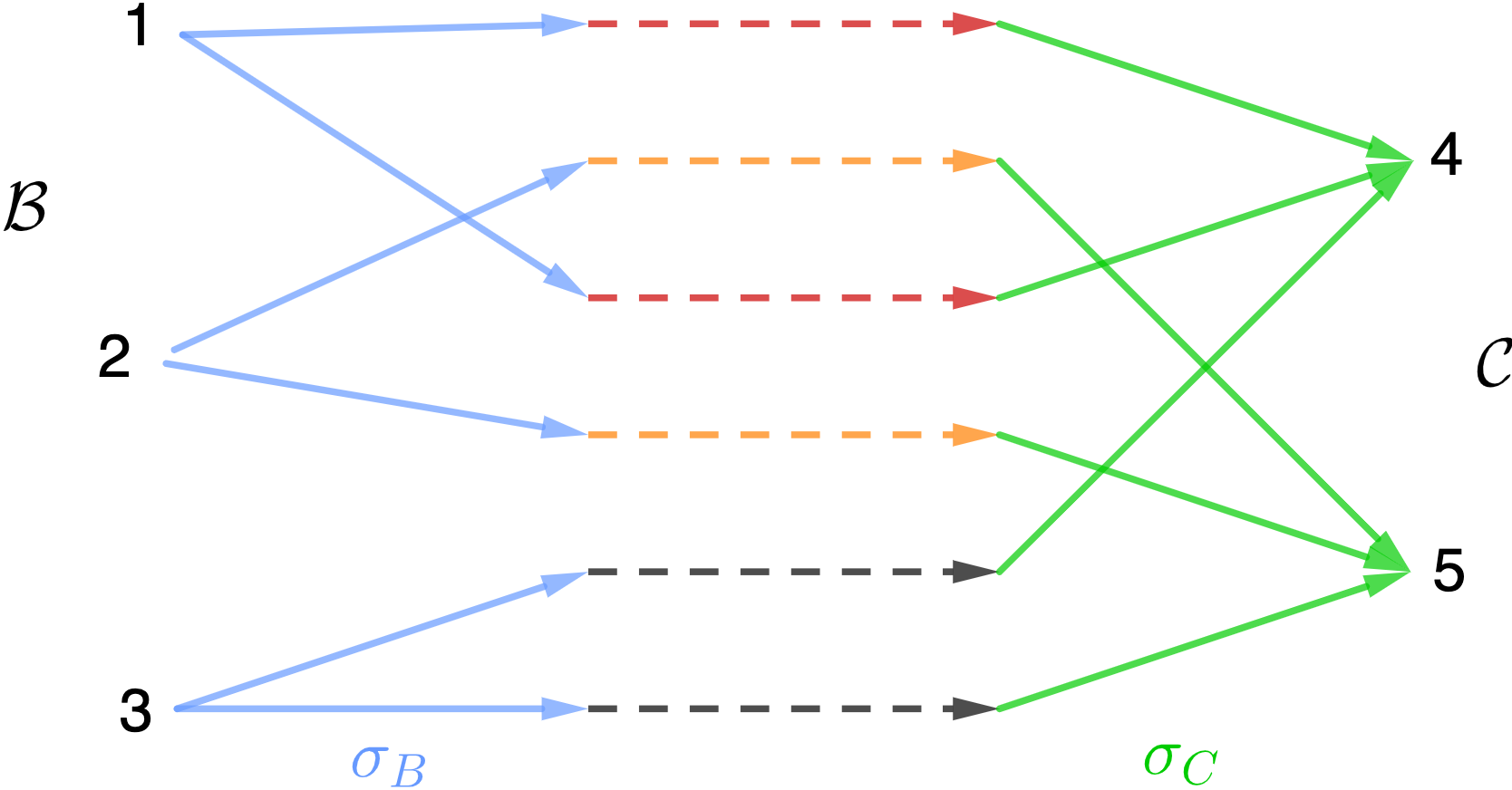}}
    \caption{Each bijective assignation $\sigma$ is composed of a $\sigma_B$ and a $\sigma_C$: see (a), but not every $(\sigma_B, \sigma_C)$ pair gives a bijective assignation $\sigma$. In (b), the orange and red edges are assigned to the same (excitation, measurement) pair: $\sigma_B$ and $\sigma_C$ are not compatible.}
    \label{fig:sigma}
\end{figure}

Besides, the determinant is expressed as the sum over all possible row-column permutations by the Leibniz formula:
\begin{align} \label{eq:leibniz}
    \det \hat K 
    = \sum_{\sigma \in S} \sgn(\sigma) \prod_{\alpha \in E^\Delta} T'_{\alpha, \sigma_B(\alpha)} T_{\sigma_C(\alpha), \alpha},
\end{align}
where each row-column permutation corresponds to a bijective assignation $\sigma:E^\Delta \rightarrow \mathcal{B} \times \mathcal{C}$, composed of (Fig. \ref{fig:sigma} (a)):
\begin{itemize}
    \item an excitation assignation $\sigma_B:E^\Delta \rightarrow \mathcal{B}$ in which $n_C$ unknown edges $\alpha$ are assigned to each excitation $b$,
    \item a measurement assignation $\sigma_C:E^\Delta \rightarrow \mathcal{C}$ in which $n_B$ unknown edges $\alpha$ are assigned to each measurement $c$.
\end{itemize}
Each bijective assignation $\sigma$ is composed of a $\sigma_B$ and a $\sigma_C$, but not every $(\sigma_B, \sigma_C)$ pair gives a bijective assignation $\sigma$, e.g. see Fig. \ref{fig:sigma} (b). We say that $\sigma_B$ and $\sigma_C$ are compatible if they form a bijection.
$S$ is the set composed of all bijective assignations $\sigma$, and $\sgn(\sigma)$ equals $+1$ if the number of transpositions in assignation $\sigma$ is even, and $-1$ otherwise. \new{A transposition is the swap of two elements. Each $\sigma$ is obtained by combining a certain number of transpositions (although such decomposition is not unique, the number of transpositions always has same parity).}

\medskip

The graph-theoretical conditions of this section rely on \emph{vertex-disjoint paths}: we say that a group of paths are vertex-disjoint if no two paths of this group contain the same vertex \cite{hendrickx2018identifiability}, see Fig. \ref{fig:vdp}. The following lemma links vertex-disjoint paths in the graph with the generic rank of transfer matrices:

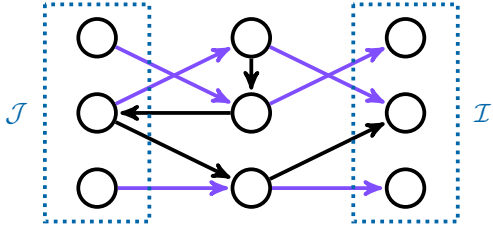
\begin{figure}
\centering
\begin{tikzpicture}[->,>=stealth',shorten >=1pt,auto,node distance=1cm,line width=0.5mm,
	  state/.style={circle, draw, minimum size=5mm}]
	  \node[state] (C)				{};
	  \node[state] (B) [above of=C] {};
	  \node[state] (D) [below of=C] {};
	  \node[state] (E) [right = 1.5 cm of B]  {};
	  \node[state] (F) [right = 1.5 cm of C] {};
	  \node[state] (G) [right = 1.5 cm of D] {};
	  \node[state] (H) [right = 1.5 cm of E]  {};
	  \node[state] (I) [right = 1.5 cm of F] {};
	  \node[state] (J) [right = 1.5 cm of G] {};
	  \path (B) edge [color = myPurple]   node {} (F)
			(F) edge [color = myPurple]   node {} (H)
			(C) edge [color = myPurple]   node {} (E)
			(E) edge [color = myPurple]   node {} (I)
			(D) edge [color = myPurple]   node {} (G)
			(G) edge [color = myPurple]   node {} (J)
			(F) edge 					  node {} (C)
			(E) edge 					  node {} (F)
			(C) edge 					  node {} (G)
			(G) edge 					  node {} (I);
   	 \draw[myMidBlue,dotted,line width=0.5mm] ($(H.north west)+(-0.5,0.25)$)  rectangle ($(J.south east)+(0.5,-0.25)$);
   	 \node (I) [right = 0.5 cm of I] {\textcolor{myMidBlue}{$\mathcal{I}$}};
   	 \draw[myMidBlue,dotted,line width=0.5mm] ($(B.north west)+(-0.5,0.25)$)  rectangle ($(D.south east)+(0.5,-0.25)$);
   	 \node (J) [left = 0.5 cm of C] {\textcolor{myMidBlue}{$\mathcal{J}$}};
\end{tikzpicture}
\caption{One can route up to three vertex-disjoint paths from $\mathcal{J}$ to $\mathcal{I}$ (in purple), i.e. $\beta_{\mathcal{J} \rightarrow \mathcal{I}}=3$.}
\label{fig:vdp}
\end{figure}

\medskip

\begin{lemma} \cite{van1991graph, hendrickx2018identifiability} \label{lemma:grank_vdp}
    Let $T_{\mathcal{I}, \mathcal{J}}$ be the transfer matrix from nodes of set $\mathcal{J}$ to nodes of set $\mathcal{I}$. We have
    \begin{align} \label{eq:rank_vdp}
        \grank T_{\mathcal{I}, \mathcal{J}} = \beta_{\mathcal{J} \rightarrow \mathcal{I}},
    \end{align}
    where \newg{the generic rank of $T_{\mathcal{I}, \mathcal{J}}$ is the rank for \emph{almost all} $G = G^0 + G^\Delta$, }
    and $\beta_{\mathcal{J} \rightarrow \mathcal{I}}$ stands for the maximum number of vertex-disjoint paths that can be routed from $\mathcal{J}$ to $\mathcal{I}$.
\end{lemma}

\medskip

\newg{Note that the rank is taken for almost all unknown transfer functions $G^\Delta$ \emph{and almost all} known transfer functions $G^0$, otherwise {\eqref{eq:rank_vdp}} would not hold for problematic values of $G^0$, as shown in Example 1 of }
 \cite{shi2020generic}.

\medskip

\newo{Combining Leibniz formula {\eqref{eq:leibniz}} with Lemma {\ref{lemma:grank_vdp}} yields the following basic proposition. An assignation $\sigma$ is \emph{connected} if for each unknown edge $\alpha$ there is a path from its assigned excitation $\sigma_B(\alpha)$ to its assigned measurement $\sigma_C(\alpha)$, in which the unknown edge $\alpha$ is included.}

\medskip

\begin{prop} \label{prop:weak_condition}
A necessary and a sufficient condition:
\begin{itemize}
    \item If a network is generically decoupled-identifiable, then there is \emph{at least} one connected bijective assignation $\sigma$.
    \item If there is \emph{only} one connected bijective assignation $\sigma$, then the network is generically decoupled-identifiable.
\end{itemize}
\end{prop}

\medskip

\newo{Proposition {\ref{prop:weak_condition}} relies on {\eqref{eq:leibniz}}, but this expression can be further developed, and algebraic manipulations allow to derive a stronger condition in terms of vertex-disjoint paths.
The following lemma is the main building brick of this stronger condition. }
 In the following lemma, an assignation $\sigma_B : E^\Delta \rightarrow \mathcal{B}$ (respectively $\sigma_C : E^\Delta \rightarrow \mathcal{C})$ is \emph{connected} if there is a path between each unknown edge $\alpha$ and its assigned excitation $\sigma_B(\alpha)$ (resp. measurement $\sigma_C(\alpha)$).

\medskip

\begin{lemma} \label{lemma:excitation_assignation}
    If a network is generically decoupled-identifiable, then there is at least one assignation $\sigma_B$ (resp. $\sigma_C$) such that:
    \begin{enumerate}[(a)]
        \item $n_C$ (resp. $n_B$) unknown edges $\alpha$ are assigned to each excitation $b$ (resp. measurement $c$)
        \item $\sigma_B$ (resp. $\sigma_C$) is connected \label{itm:connected}
        \item for each excitation $b$ (resp. measurement $c$), there are $n_C$ (resp. $n_B$) vertex-disjoint paths between the edges assigned to $b$ (resp. $c$) and measurements $\mathcal{C}$ (resp. excitations $\mathcal{B}$). \label{itm:vdp}
    \end{enumerate}
    If there is only one such assignation, then this condition is also sufficient for generic decoupled identifiability.
\end{lemma}

\medskip

Combining Lemma \ref{lemma:excitation_assignation} for $\sigma_B$ and $\sigma_C$ yields the main result of this paper.

\medskip

\begin{thm} \label{thm:strong_condition}
    If a network is generically decoupled-identifiable, then there is at least one assignation $\sigma$ such that:
    \begin{enumerate}[(a)]
        \item \label{itm:thm_sigmaB} $n_C$ unknown edges $\alpha$ are assigned to each excitation $b$
        \item $n_B$ unknown edges $\alpha$ are assigned to each measure $c$ \label{itm:thm_sigmaC}
        \item $\sigma$ is connected
        \item for each excitation $b$, there are $n_C$ vertex-disjoint paths between the edges assigned to $b$ and measurements $\mathcal{C}$. \label{itm:thm_vdp_C}
        \item for each measurement $c$, there are $n_B$ vertex-disjoint paths between edges assigned to $c$ and excitations $\mathcal{B}$. \label{itm:thm_vdp_B}
    \end{enumerate}
    If there is only one such assignation, then this condition is also sufficient for generic decoupled identifiability.
\end{thm}

%% file: 5-dicussion.tex
First, we highlight some subtle difference between the conditions of Theorem \ref{thm:strong_condition} and those of Proposition \ref{prop:weak_condition}:
\begin{enumerate}[(i)]
    \item In contrast to Proposition \ref{prop:weak_condition}, the assignation $\sigma$ of Theorem \ref{thm:strong_condition} is not necessarily bijective: two edges assigned to the same excitation can be assigned to the same measurement.
    \item The assignations of condition \ref{itm:thm_sigmaB} do not necessarily match the vertex-disjoint paths of condition \ref{itm:thm_vdp_B}.
    \item The assignations of condition \ref{itm:thm_sigmaC} do not necessarily match the vertex-disjoint paths of condition \ref{itm:thm_vdp_C}.
\end{enumerate}

We believe that there could be a stronger version of Theorem \ref{thm:strong_condition}, in which $\sigma$ is bijective, as in Proposition \ref{prop:weak_condition}. It would unify Proposition \ref{prop:weak_condition} with Theorem \ref{thm:strong_condition}, and extend the vertex-disjoint path conditions of \cite{hendrickx2018identifiability}.

\medskip

Besides, Proposition \ref{prop:weak_condition} and Theorem \ref{thm:strong_condition} provide necessary conditions for generic decoupled identifiability and sufficient ones.
As shown in Proposition \ref{prop:bili_necessary}, generic decoupled identifiability is necessary for generic local identifiability, which was itself shown to be necessary for generic identifiability in \cite{legat2020local}. 
Therefore, the necessary conditions derived in Section \ref{sec:conditions} apply to (generic) (local) identifiability.

\medskip

Regarding sufficiency, we remind that in the systematic tests we have conducted (code at \cite{matlab}), we have found no network which is generically decoupled-identifiable but not generically locally identifiable. No counterexample to the sufficiency of generic local identifiability for generic identifiability has been found either yet \cite{legat2020local}, but we could not conduct numerical tests for this one since we lack of a general criterion to check generic identifiability.
Consequently, one might hope that the sufficient conditions of Section \ref{sec:conditions} apply to generic (local) identifiability, but this remains an open question.

%% file: 6-conclusion.tex
This work was motivated by one main open question: determining path-based conditions for generic local identifiability in networked systems.

\medskip

The decoupled version of generic identifiability we have introduced allowed to look at generic identifiability from a new angle, on a larger graph which decouples excitations and measurements. In particular, we have derived necessary conditions in terms of vertex-disjoint paths in the larger graph and sufficient ones. The necessary conditions extend to generic (local) identifiability, but whether the sufficient conditions extend as well remains an open question.

\medskip

We believe that the two identifiability conditions of this paper could be merged into a unifying stronger condition. It would extend results of full excitation or measurement. 

\medskip

A further open question would be to obtain graph theoretical conditions for local identifiability of a subset of edges when not all edges can be recovered, a question for which \cite{legat2020local} gives an algebraic necessary and sufficient condition.

%% file: 7-appendix.tex
\subsection{Proofs of Section \ref{sec:decoup-identif}.}

\medskip

\textbf{Proof of Proposition \ref{prop:rank_hat_K_bilidentif}:}
    We first prove that \eqref{eq:CTDeltaT'B_bili} holds if and only if $\hat K(G,G')$ has full rank:
    vectorizing \eqref{eq:CTDeltaT'B_bili} gives
    \begin{align*}
        \prt{B^\top T(G')^\top \otimes C \, T(G)} \myvec(\Delta) = 0
        \Rightarrow \myvec(\Delta) = 0,
    \end{align*}
    where $\myvec$ stacks the columns of its $|\mathcal{B}|\times|\mathcal{C}|$ matrix argument into a $|\mathcal{B}|\cdot|\mathcal{C}|$ vector.
	Then, observing that $\myvec \prt{\Delta} = I_{G^\Delta} \, \delta$ yields $ \hat K \delta = 0 \Rightarrow \delta = 0 \quad \forall \, \delta \in \mathbb{C}^{|E^\Delta|}$, which gives that $\hat K(G,G')$ has full column rank $|E^\Delta|$.

	Then, we prove that the full-rankness of $\hat K(G,G')$ is a generic property. The following argument is inspired from the proof of Lemma B.1 in \cite{legat2020local}:
	the set on which the rank of $\hat K(G,G')$ is not full can be expressed as the intersection of zero sets of determinants of submatrices {\cite{van1991graph}}, which are analytic in their entries $(G,G')$. Since nonconstant analytic functions vanish only on a lower-dimensional set {\cite{d2010introduction}}, the set on which $\hat K(G,G')$ does not reach its maximal rank has lower dimension.
	This extends the notion of generic rank introduced in {\cite{davison1977connectability}}.
\hfill \rule{2mm}{2mm}

\bigskip

\textbf{Proof of Proposition \ref{prop:bili_necessary}:}
    \newg{Consider a generically locally identifiable network. Take $G$, at which the network is locally identifiable. }
     Equation \eqref{eq:CTDeltaTB_Delta_net} holds, and it is decoupled-identifiable at $(G,G'=G)$ by definition.
    \newg{Proposition {\ref{prop:rank_hat_K_bilidentif}} asserts that decoupled identifiability is a generic property, thus the fact it holds for one variable $(G,G'=G)$ implies it holds for almost all $(G,G')$.}
\hfill \rule{2mm}{2mm}

\bigskip

\textbf{Proof of Proposition \ref{prop:equivalence_bili_decoupled}:}
    \newo{From definitions {\ref{def:global_identif}} and {\ref{def:decoupled_network}}, the decoupled network $\hat G(G,G')$ is generically identifiable if, for all $\tilde{G}^\Delta$ with same zero entries as $G^\Delta$, there holds}
    \begin{align} \label{eq:dec-net-identif}
		\hat C \, (I - \hat{\tilde{G}})^{-1} \, \hat B = \hat C \, (I-\hat G)^{-1} \, \hat B \Rightarrow \tilde{G}^\Delta = G^\Delta
	\end{align}
	\newo{for almost all $G^\Delta$ with same zeros, where }
	$
        \hat{\tilde{G}} =
        \begin{bmatrix}
            G & \tilde G^\Delta \\
            0 & G'
        \end{bmatrix}.
    $
    \newo{Developing {\eqref{eq:dec-net-identif}} yields}
    \begin{align*}
        C \, T(G) \, \tilde{G}^\Delta \, T(G') \, B =
        C \, T(G) \, G^\Delta \, T(G') \, B \Rightarrow
        \tilde{G}^\Delta = G^\Delta,
    \end{align*}
    \newo{and bringing out common terms gives}
    \begin{align*}
        C \, T(G) \, \underbrace{(\tilde{G}^\Delta - G^\Delta)}_{\triangleq \Delta} \, T(G') \, B = 0 \Rightarrow
        \underbrace{\tilde{G}^\Delta - G^\Delta}_{= \Delta} = 0,
    \end{align*}
    \newo{which matches Definition {\ref{def:decoupled_identif}}, and requiring this to hold for almost all $(G,G')$ completes the proof.}
\hfill \rule{2mm}{2mm}

\medskip

\subsection{Proofs of Section \ref{sec:conditions}}

\medskip

\textbf{Proof of Proposition \ref{prop:weak_condition}:}
    Consider a generically decoupled-identifiable network. Corollary \ref{corol:det_hat_K_bilidentif} yields that $\det \hat K \neq 0$ generically.
    Besides, the Leibniz formula for the determinant gives \eqref{eq:leibniz}:
    \begin{align*}
        \det \hat K
        = \sum_{\sigma \in S} \sgn(\sigma) \underbrace{\prod_{\alpha \in E^\Delta} T'_{\alpha, \sigma_B(\alpha)} T_{\sigma_C(\alpha), \alpha}}_{\triangleq \tau(\sigma)}.
    \end{align*}
    
    \vspace*{-4mm}
    
    Since $\det \hat K \neq 0$ generically, there is at least one assignation $\sigma$ with a nonzero term $\tau(\sigma)$. For this assignation, each factor of the product must be generically nonzero, i.e. $T'_{\alpha, \sigma_B(\alpha)} \neq 0$ and $T_{\sigma_C(\alpha),\alpha} \neq 0$ for all unknown edge $\alpha$.
    Lemma \ref{lemma:grank_vdp} for $\mathcal{J} = \{j\}, \mathcal{I} = \{i\}$ ensures that $T_{i,j}$ is generically nonzero if and only if there is a path from $j$ to $i$.
    Hence, we can construct a path going from $\sigma_B(\alpha)$ to $\sigma_C(\alpha)$, including $\alpha$.
    
    If more than one assignation $\sigma$ have a nonzero term $\tau(\sigma)$, the condition is not sufficient since terms could cancel each other. But if there is only one such assignation, then the condition is also sufficient.
\hfill \rule{2mm}{2mm}

\medskip

The proof of Lemma {\ref{lemma:excitation_assignation}} relies on Lemma \ref{lemma:signature}, which handles the signature decomposition of the assignation.

\medskip

\begin{lemmaappendix} \label{lemma:signature}
    \new{Let $\sigma$ be a bijective assignation comprised of compatible $\sigma_B$ and $\sigma_C$. Its signature decomposes into}
    \begin{align*}
        \sgn(\sigma) = \sgn(\sigma_B) \ \sgn(\sigma_C),
    \end{align*}
    \new{where $\sgn(\sigma_B) \triangleq \prod_c \sgn(\sigma_{B,c}), \ \sgn(\sigma_C) \triangleq \prod_b \sgn(\sigma_{C,b}),$ $\sigma_{C,b}$ is the sub-assignation of $\sigma_C$ for the edges $\alpha$ assigned to the same excitation $b$ in $\sigma_B$ (see Fig. {\ref{fig:subsigmaC}} and equation {\eqref{eq:def:sigmaCb}}), and $\sigma_{B,c}$ is defined analogously.}
\end{lemmaappendix}

\medskip

\textbf{Proof:}
    \new{An assignation $\sigma$ is comprised of a sequence of transpositions that swap two edges: denote them $\alpha_{1,1}$ and $\alpha_{2,2}$, initially assigned to $(b_1,c_1)$ and $(b_2, c_2)$ respectively. Since $\sigma$ is bijective, either (i) $b_1 = b_2, c_1 \neq c_2$, (ii) $b_1 \neq b_2, c_1 = c_2$, or (iii) $b_1 \neq b_2, c_1 \neq c_2$.
    Besides, every transposition of case (iii) can be decomposed into transpositions of type (i) and (ii), as shown in Table {\ref{tab:decomposition_transposition}}.}
    \begin{table}[H]
        \centering
        \begin{tabular}{c|c c c c c c c}
            $(b_1, c_1)$ & $\boxed{\alpha_{1,1}}$ & $\xrightarrow[]{\text{(i)}}$ & $\alpha_{1,2}$ & $\xrightarrow[]{\text{(ii)}}$ & $\boxed{\alpha_{1,2}}$ & $\xrightarrow[]{\text{(i)}}$ & $\alpha_{2,2}$ \\
            & $\updownarrow$ & & & & $\updownarrow$ & & \\
            $(b_1, c_2)$ & $\boxed{\alpha_{1,2}}$ & & $\boxed{\alpha_{1,1}}$ & & $\boxed{\alpha_{2,2}}$ & & $\alpha_{1,2}$ \\
            $(b_2, c_1)$ & $\alpha_{2,1}$ &  & $\alpha_{2,1} \ \updownarrow$ &  & $\alpha_{2,1}$ &  & $\alpha_{2,1}$ \\
            $(b_2, c_2)$ & $\alpha_{2,2}$ && $\boxed{\alpha_{2,2}}$ && $\alpha_{1,1}$ && $\alpha_{1,1}$
        \end{tabular}
        \caption{\new{Decomposition of a class-(iii) transposition into class-(i) and -(ii) transpositions. The boxes and arrows indicate the swap to the next column.}}
        \label{tab:decomposition_transposition}
    \end{table}
    
    \vspace*{-2mm}
    
    \new{Hence, any assignation $\sigma$ can be decomposed into a sequence of class-(i) and -(ii) transpositions only. Take one such sequence, and denote its number of transpositions $\pi(\sigma)$.
    Within this sequence, consider the transpositions of class (i) with same fixed $b$. Combining these transpositions gives the sub-assignation $\sigma_{C,b}$, defined in equation {\eqref{eq:def:sigmaCb}} and illustrated in Fig. {\ref{fig:subsigmaC}}. Hence, the number of such transpositions gives the parity of $\sigma_{C,b}$, and we denote it by $\pi(\sigma_{C,b})$. The same holds for all $b$, and for the class-(ii) transpositions, for each $c$. Enumerating all $\pi(\sigma)$ transpositions of the considered sequence allows to rewrite it as a sum:}
    \begin{align*}
        \pi(\sigma) = \sum_c \pi(\sigma_{B,c}) + \sum_b \pi(\sigma_{C,b}).
    \end{align*}
    \new{Since $\sgn(\sigma) \triangleq (-1)^{\pi(\sigma)}$, this yields
    }
    \begin{align*}
        \qquad \quad \
        \sgn(\sigma) = \prod_c (-1)^{\pi(\sigma_{B,c})}
        \ \prod_b (-1)^{\pi(\sigma_{C,b})}.
        \qquad \quad \
        \rule{2mm}{2mm}
    \end{align*}
    
\newpage

\textbf{Proof of Lemma \ref{lemma:excitation_assignation}:}
    This proof is derived for $\sigma_B$, but an analogous argument can be made for $\sigma_C$.
    Consider a generically decoupled-identifiable network. Corollary \ref{corol:det_hat_K_bilidentif} yields that $\det \hat K \neq 0$ generically.
    Besides, the Leibniz formula for the determinant gives \eqref{eq:leibniz}:
    \begin{align*}
        \det \hat K 
        &= \sum_{\sigma \in S} \sgn(\sigma) \prod_{\alpha \in E^\Delta} T'_{\alpha, \sigma_B(\alpha)} T_{\sigma_C(\alpha), \alpha}\\
        &\, \downarrow \quad \text{grouping $\sigma$ with same $\sigma_B$, \new{Lemma {\ref{lemma:signature}} for }}
        \sgn \\
        &= \sum_{\sigma_B} \sgn(\sigma_B) \sum_{\sigma_C \in S_C^{\sigma_B}} \sgn(\sigma_C) \prod_{\alpha} T'_{\alpha, \sigma_B(\alpha)} T_{\sigma_C(\alpha) \alpha},
    \end{align*}
    where $S_C^{\sigma_B}$ is the set of all $\sigma_C$ compatible with $\sigma_B$, i.e. such that $(\sigma_B,\sigma_C)$ is bijective: if two edges are assigned to the same $b$ by $\sigma_B$, their $c$ assigned by $\sigma_C$ must be different.
    
    \medskip
    
    Bringing out common $T'$ terms yields $\det \hat K =$
    \begin{align*}
        \sum_{\sigma_B} \sgn(\sigma_B) \prod_\alpha T'_{\alpha, \sigma_B(\alpha)}
        \underbrace{\sum_{\sigma_C \in S_C^{\sigma_B}} \sgn(\sigma_C) \prod_{\alpha} T_{\sigma_C(\alpha), \alpha}}_{\triangleq \gamma_{\sigma_B}}
    \end{align*}
    We develop further $\gamma_{\sigma_B}$, the factor accounting for the assignations $\sigma_C$ compatible with $\sigma_B$. We decompose the product on all unknown edges $\alpha$ into subgroups of $\alpha$ assigned to the same excitation $b$ in $\sigma_B$ (subgroups denoted as $\sigma_B^{-1}(b)$):
    \begin{align} \label{eq:before_key_step}
        \gamma_{\sigma_B}
        = \sum_{\sigma_C \in S_C^{\sigma_B}} \sgn(\sigma_C) \prod_{b} \prod_{\alpha \in \sigma_B^{-1}(b)} T_{\sigma_C(\alpha), \alpha}
    \end{align}
    The assignation $\sigma_C$ can be decomposed into $n_B$ sub-assignations for the edges $\alpha$ assigned to the same excitation $b$ in $\sigma_B$ (i.e. the $\alpha$ belonging to $\sigma_B^{-1}(b)$, see Fig. \ref{fig:subsigmaC}): 
    
    \begin{figure}
        \centering
        \includegraphics[width=\linewidth]{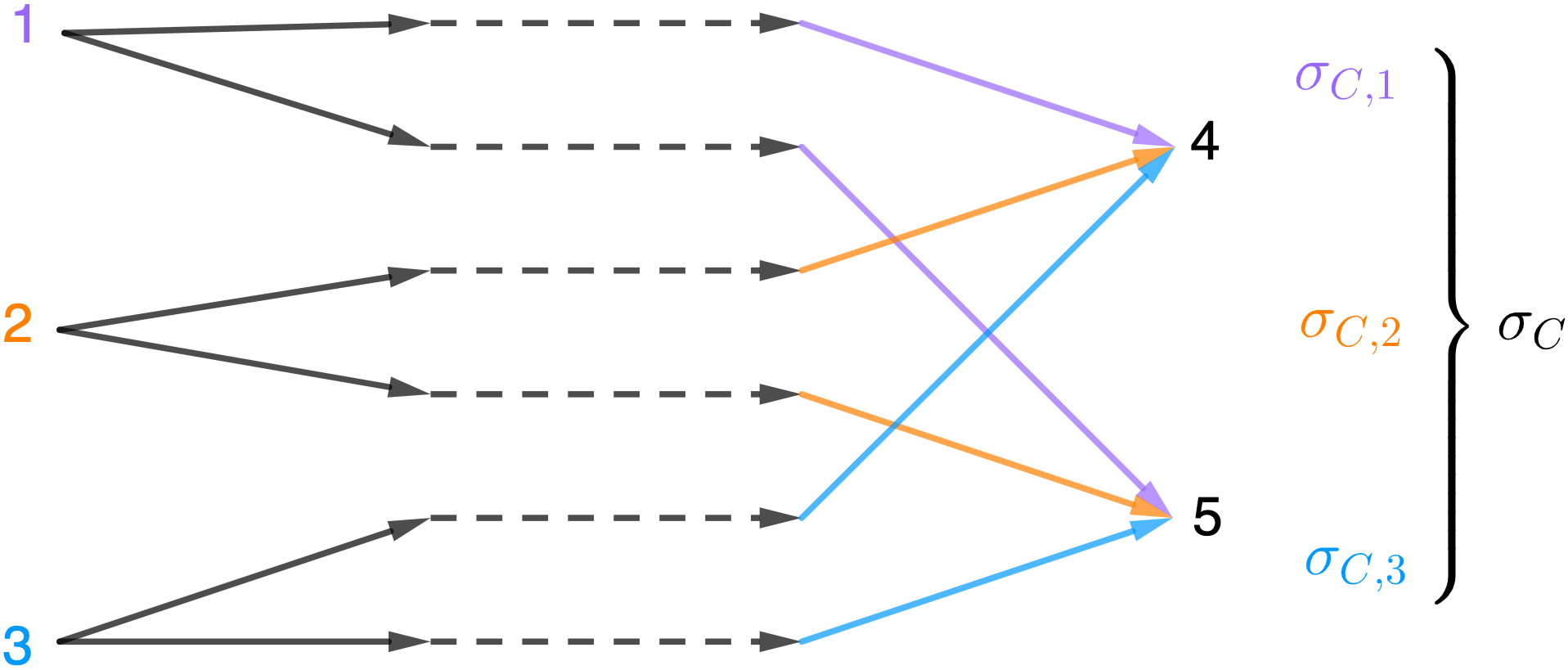}
        \caption{The assignation $\sigma_C$ can be decomposed into $n_B$ sub-assignations for the edges $\alpha$ assigned to the same excitation $b$ in $\sigma_B$.}
        \label{fig:subsigmaC}
    \end{figure}
    
    \begin{align} \label{eq:def:sigmaCb}
        \sigma_{C,b}: \sigma_B^{-1}(b) \rightarrow \mathcal{C}: \alpha \rightarrow \sigma_{C,b}(\alpha)
    \end{align}
    Here is the key step of the proof: since those sub-assignations $\sigma_{C,b}$ do not depend on each other, \eqref{eq:before_key_step} can be rewritten by regrouping in the same factor the edges $\alpha$ assigned to the same excitation $b$ in the assignation $\sigma_B$:
    \begin{align} \label{eq:after_key_step}
        \gamma_{\sigma_B}
        = \prod_{b} \sum_{\sigma_{C,b} \in S_{C,b}^{\sigma_B}} \sgn(\sigma_{C,b}) \prod_{\alpha \in \sigma_B^{-1}(b)} T_{\sigma_{C,b}(\alpha), \alpha},
    \end{align}
    where $S^{\sigma_B}_{C,b}$ denotes the set of sub-assignations $\sigma_{C,b}$ compatible with $\sigma_B$.
    
    Furthermore, each sub-assignation $\sigma_{C,b}$ compatible with $\sigma_B$ is bijective: 
    since each pair $(b,c)$ must be assigned to only one edge $\alpha$ for compatibility, each $\alpha$ of $\sigma_B^{-1}(b)$ must be assigned to a different measurement, i.e. $\sigma_{C,b}$ must be injective. Besides, $\sigma_B^{-1}(b)$ is composed of $n_C$ edges, so $\sigma_{C,b}$ is bijective.
        
    \medskip
    
    All compatible $\sigma_{C,b}$ must be taken in order to cover all assignations $\sigma$.
    Thus $S^{\sigma_B}_{C,b}$ contains all possible bijections $\sigma_{C,b}$.
    Hence, each factor of \eqref{eq:after_key_step} is a determinant:
    \begin{align*}
        \gamma_{\sigma_B}
        = \prod_{b} \det T_{\mathcal{C},  \sigma_B^{-1}(b)}.
    \end{align*}
    Putting it all together, we have
    \begin{align} \label{eq:dethatK_detT}
        \det \hat K
        = \sum_{\sigma_B} \sgn(\sigma_B) \underbrace{\prod_\alpha T'_{\alpha, \sigma_B(\alpha)}
        \prod_{b} \det T_{\mathcal{C}, \sigma_B^{-1}(b)}}_{\triangleq \tau_B(\sigma_B)}
    \end{align}
    Since $\det \hat K \neq 0$ generically, there is at least one assignation $\sigma_B$ with a nonzero term $\tau_B(\sigma_B)$. For this assignation, each factor of the product must be nonzero, i.e. (i) $T'_{\alpha, \sigma_B(\alpha)} \neq 0$ for all unknown edge $\alpha$ and (ii) $\det T_{\mathcal{C}, \sigma_B^{-1}(b)} \neq 0$ for all excitation $b$.
        
    \medskip
    
    Lemma \ref{lemma:grank_vdp} for single nodes $\mathcal{J} = \{j\}, \mathcal{I} = \{i\}$ ensures that $T_{i,j}$ is generically nonzero if and only if there is a path from $j$ to $i$, so (i) gives condition \ref{itm:connected}.
        
    \medskip
    
    Besides, (ii) is equivalent to $\rank T_{\mathcal{C}, \sigma_B^{-1}(b)} = n_C$ for all excitation $b$, and Lemma \ref{lemma:grank_vdp} yields condition \ref{itm:vdp}.
    
    \medskip
    
    If more than one assignation $\sigma$ have a nonzero term $\tau_B(\sigma_B)$, the condition is not sufficient since terms could cancel each other. But if there is only one such assignation, the condition is also sufficient.
\hfill \rule{2mm}{2mm}

\bigskip

\textbf{Proof of Theorem \ref{thm:strong_condition}:}
    Consider a generically decoupled-identifiable network. Corollary \ref{corol:det_hat_K_bilidentif} yields that $\det \hat K \neq 0$ generically.
    Combine \eqref{eq:dethatK_detT} with its analogous version for the measurements:
    \begin{align*}
        \prt{\det \hat K}^2
        &= \prt{\sum_{\sigma_B} \sgn(\sigma_B) \ \tau_B(\sigma_B)}
        \prt{\sum_{\sigma_C} \sgn(\sigma_C) \ \tau_C(\sigma_C)}\\
        &= \sum_{\sigma_B} \sum_{\sigma_C} \sgn(\sigma_B, \sigma_C)\ \tau_B(\sigma_B)\ \tau_C(\sigma_C),
    \end{align*}
    where $\tau_B(\sigma_B)$ stands for the contribution of assignation $\sigma_B$ in $\det \hat K$ as defined in \eqref{eq:dethatK_detT}, and $\tau_C(\sigma_C)$ is defined analogously.
        
    \medskip
    
    Since $\det \hat K$ is generically nonzero, so is $(\det \hat K)^2$, and there is at least one pair of assignations $(\sigma_B, \sigma_C)$ such that $\tau_B(\sigma_B) \neq 0$ and $\tau_C(\sigma_C) \neq 0$.
    Hence, the connectivity and vertex-disjoint path conditions follow from Lemma \ref{lemma:excitation_assignation}.
    
    \medskip
    
    If more than one pair $(\sigma_B, \sigma_C)$ have a nonzero term $\tau_B(\sigma_B) \, \tau_C(\sigma_C)$, the condition is not sufficient since terms could cancel each other. But if there is only one such pair, the condition is also sufficient.
\hfill \rule{2mm}{2mm}